\documentclass[12pt,a4paper]{article}
\usepackage{amsmath,amstext,amssymb,amscd,tikz-cd, color}
\usepackage{graphicx}

\usepackage[english]{babel}

\newcommand{\Xcomment}[1]{}

\newtheorem{claim}{Claim}
\newtheorem{thm}{Theorem}
\newtheorem{lemma}{Lemma}
\newtheorem{cor}{Corollary}

\newtheorem{theorem}{Theorem}

\makeatletter \@addtoreset{equation}{section} \makeatother

{\hfill$\qed$\medskip}

\def\qed{ \ \vrule width.1cm height.3cm depth0cm}

 \makeatletter
\renewcommand{\section}{\@startsection{section}{1}{0pt}%
{-3.5ex plus -1ex minus -.2ex}{2.3ex plus .2ex}%
{\normalfont\Large}}
 \makeatother

 \makeatletter
\renewcommand{\subsection}{\@startsection{subsection}{2}{0pt}%
{-3.0ex plus -1ex minus -.2ex}{-1.5ex plus .2ex}%
{\normalfont\large\bf}}
 \makeatother

 \Xcomment{
 \makeatletter
\renewcommand{\subsection}{\@startsection{subsection}{2}{0pt}%
{-3.0ex plus -1ex minus -.2ex}{1.5ex plus .2ex}%
{\normalfont\large\bf}}
  }
 \makeatother
 
\newcommand{\iso}{\overset{\sim}{\longrightarrow}}
\newcommand{\isom}{\overset{\sim}{=}}
\newcommand{\lra}{\longrightarrow}

\newcommand{\lla}{\longleftarrow} 

\def\BC{{\mathbb{C}}}

\def\BP{{\mathbb{P}}}

\def\R{{\mathbb{R}}}

\def\BQ{{\mathbb{Q}}}
\def\BR{{\mathbb{R}}}

\def\BZ{{\mathbb{Z}}}

\def\CM{{\cal M}}

\def\CP{{\cal P}}

\def\fC{\mathfrak{C}}

\def\fD{\mathfrak{D}}

\def\tilde{\widetilde}

\def\bar{\overline}

\def\tC{\tilde{C}}
\def\tD{\tilde{D}}
\def\tE{\tilde{E}}
\def\tp{\tilde{p}}

\newcommand{\be}{\begin{equation}}
\newcommand{\ee}{\end{equation}}

\begin{document}

\title{PENTAGRAMMA MIRIFICUM. II}

\author{Vadim Schechtman}

\date{\today}

\maketitle

\section{Introduction}

\

\subsection{Main topics} This article is a continuation of \cite{S}. We discuss here the 
following topics. 

(a) The moduli space of miraculous pentagrams as the del Pezzo surface $S_5$ of degree $5$.
 
(b) Relation to the Vinberg's most algebraic $K3$ surface, \cite{V}.
 
(c) Some particular miraculous pentagrams which we call 
{\it Fibonacci pentagrams}.

(d) Fermat's ascent and the Jacobi's proof of the  Poncelet theorem. 

\subsection{The surface $S_5$ and the root system $A_4$} The del Pezzo surface $S_5$ of degree $5$ may be defined as the blowup of 
$\BP^2$ at $4$ points, say $p_1, \ldots, p_4$, in general position, cf. \cite{M}, Ch. IV, \S 2, Thm. 2.5.

 It contains $10$ lines: $6$ strict transforms of the lines passing through $p_i, p_j, i\neq j$, together with 
$4$ exceptional divisors; denote these lines $\ell_i, 1\leq i\leq 10$. They meet at $15$ points; their 
incidence graph is the Petersen graph $P$, \cite{P}. 


In the language of \cite{M}, Ch. IV the surface $S_5$ corresponds to the root system $R = R_4$ of type $A_4$, so 
$$ 
R\subset H^2(S_5; \Bbb Z) \subset H^2(S_5; \BR)\isom \BR^5.
$$
(we identify $Pic(S_5)$ with $H^2(S_5; \Bbb Z)$ using the first Chern class).

The Weil group $W(R)$ is the 
symmetric group $\Sigma_5$ which coincides with $Aut(S_5) = Aut(P)$. 

In \cite{M}, Def. 1.7 and Thm. 1.8 the root system $R$ and the set of exceptional curves $I_4$ 
are identified explicitly. 

\subsubsection{Exercise} 

Describe the set of exceprional curves 
$$
I_4 = \{\ell_1, \ldots, \ell_{10}\}\subset H^2(S_5; \Bbb Z)\isom \Bbb Z^5 
$$
explicitly using Manin's description in \cite{M}, Ch. IV, \S 1, Def. 1.7. See the Petersen graph therein.

{\bf Solution.} See \cite{M}, Ch. IV, \S 4, Fig. 4. The Petersen graph is $\Gamma_4$, in a somewhat unusual shape.

\subsection{Three modular interpretations}   

The surface $S_5$ admits three modular incarnations:

\ 


(a) it is 
a compatification of the space $\CM_{0,5}$ of $5$-tuples of points on $\BP^1$;

\ 

(b)  it is the moduli space $\CP_5$ of {\it miraculous pentagrams} - certain $5$-tuples of points on 
$\BP^2$.

\

The isomorphism between these two models is provided by the Veronese embedding
$$
v:\ \BP^1 \lra \BP^2,\ v(x:y) = (x^2:xy:y^2).
$$

\

(c)  The surface $\CP_5$ is the complexification of a surface $\BR\CP_5$ over $\BR$.  
According to Gauss \cite{Ga,S}, $\BR\CP_5$ may be interpreted as the moduli space of {\it miraculous  pentagrams} --- certain
$5$-tuples of points on the sphere $S^2$. Such pentagram is defined by $5$ numbers 
$\alpha, \beta, \gamma, \delta, \epsilon$ satisfying $5$ simple equations, see \eqref{eq:frieze} in 2.1  below. 
These equations define what is now called the {\it cluster algebra of type $A_2$}, cf. \cite[Example 1.3]{FZ}.

For a point $p\in \CP_5$
we call the numbers $\alpha, \beta, \gamma, \delta, \epsilon$  the {\it Gauss coordinates} 
of $p$.

\

Jacobi has associated to a point $p\in \CP_5$ an elliptic curve $E(p)$; its explicit equation in the Legendre form was written down by Gauss. On the other hand we can see the same curve in (a) realization as well.

\subsection{Vinberg $K3$ surface and the Poncelet-Jacobi pencil} 
There is a remarkable $K3$ surface $X_4$ introduced by Vinberg \cite{V}; it can be defined as follows. 
Take four points in general position $p_1, \ldots, p_4\in \BP^2$. We have six straight lines in $\BP^2$ 
connecting the pairs of distinct points $p_i, p_j$. Let $P\lra \BP^2$ be the double covering 
branched at these six lines; $P$ is singular at the inverse images of $p_i$. If we blowup 
them we get a nonsingular surface which is $X_4$, see \ref{vinberg-def} below.

Alternatively, we can consider the double covering of $S_5$ ramified at ten lines $\ell_i$, and then blow up their $15$ points of intersection; the resulting surface will be isomorphic to $S_4$.

Thus we have the projection  
$$
\pi: X_4 \lra S_5
$$
For the analog of the Peterson graph in $X_4$ see \ref{vinbergs-definition} below.

\subsubsection{Elliptic pencil in $X_4$} The del Pezzo $S_5$ admits a pencil of conics
$$
\phi': S_5 \lra \BP^1,
$$
cf. \ref{conics-on-delpezzo} below. Let
$$
\phi = \phi'\pi:\ \ X_4\lra \BP^1.
$$ 
For each $p\in S_5$ we can identify the corresponding Jacobi elliptic curve $E(p)$ with the fiber  
$$
E(p) =  \phi^{-1}\phi'(p)\subset X_4.
$$
These curves form an elliptic pencil in $X_4$,  
the {\it Poncelet pencil}, or {\it pope}, with $3$ degenerate fibers whose Kodaira type is the affine $D_6$.

\subsection{Fibonacci pentagrams} 
Among the miraculous pentagrams there is a distinguished one 
$p_0$ which corresponds to the regular pentagon. The corresponding elliptic curve $E(p_0)$ is degenerate, 
its Legendre equation is 
$$
y^2 = 1 - x^2.
$$
The Gauss coordinates of $p_0$ are all equal to the golden ratio 
$$
\alpha_0 = \frac{1 + \sqrt{5}}{2}.
$$
It admits the well-known expression as (the simplest among all existant) continuous fraction, whose 
convergeants are ratios of Fibonacci numbers.

In Section 5 we define a sequence of points $p_n\in \CP_5(\BQ)$ such that 
$$
\lim_{n\lra\infty} p_n = p_0.
$$
The Gauss coordinates of $p_n$ admit very simple rational expressions through the  Fibonacci numbers, see \eqref{eq:fib-coordinates}.

The elliptic curve 
$E(p_n)$ is defined over a real quadratic extension $\BQ_n$ of $\BQ$.

\subsection{Fermat and Poncelet} In Sections 6 - 12 we discuss 
 the Fermat's ascent method for solving the Diophantus "double equations" and its relation to the Poncelet theorem.
 
\subsection{Notation} 
We will use the notation from \cite{S}. 
Our ground field will be $\BC$. So $\BP^2$ will denote the projective plane over $\BC$; if necessary, 
we will consider $\BR\BP^2\subset \BP^2$, etc. 


\subsection{Acknowledgement} I am grateful to Alexander Kuznetsov for numerous consultations. In particular the idea to realize the Jacobi-Poncelet elliptic curves inside the Vinberg wonderful $K3$ surface belongs to him. I am grateful to Gil Bor for the collaboration 
at the initial stage of this work. 

\newpage

\section{Miraculous pentagrams and the del Pezzo surface $S_5$}

\

{\it Pasquale del Pezzo}, Duke of Caianello and Marquis of Campodisola (1859 –-- 1936). 

\subsection{A surface related to the Napier laws of spherical trigonometry}

 The laws of spherical trigonometry were  expressed by John Napier in 1614 as follows. Given a spherical right triangle,  one labels its 5 `parts'  (2 angles and 3 side lengths) by 
 $\alpha_1,\alpha_2',\alpha_3,\alpha_4',\alpha_5$, where $\alpha':=\pi/2-\alpha$.
 
Define 
$$
x_i:=\tan^2 \alpha_i, \quad i=1,\ldots, 5.
$$
Then one finds that
\be\label{eq:frieze}
\begin{gathered}
1+x_1=x_3 x_4,\quad
1+x_2=x_4 x_5,\quad
1+x_3=x_5 x_1,\\
1+x_4=x_1 x_2,\quad
1+x_5=x_2 x_3.
\end{gathered}
\ee

Below we will use the notation 
$$
(\alpha, \beta, \gamma, \delta, \epsilon) = (x_1, x_2, x_3, x_4, x_5)
$$
as well.

Note that these 5 formulas are obtained from any one of them by cyclic permutation of the five variables. 

It follows that the same formulas hold for the four other triangles obtained from the original one by cyclically permuting the five parts of the original triangle. The five triangles form a {\em pentagramma mirificum} (''miraculous pentagram''), a spherical 
right-angled five-pointed star. The hypotenuses form a self-polar pentagon (each vertex is perpendicular to the opposite edge). 

\begin{thm}\label{thm:frieze}
The  space of solutions to equations \eqref{eq:frieze}  is a smooth surface  in $\R^5$. Its closure $S\subset \BR\BP^5$ is a smooth compact surface, whose complexification  is a smooth complex compact surface  $S_\BC\subset\BC\BP^5$.  
\end{thm}

We leave the proof to the reader.

\

\subsection{Three definitions of $S_5$}

A del Pezzo surface is a smooth  projective surface $V$ with the very ample 
anticanonical bundle $\omega = \Omega_V^{-1}$. The degree of $V$ is $d = (\omega_V, \omega_V)$, cf. \cite{M}. 

We will be interested in the del Pezzo surface $S_5$ of degree $5$. It is unique up to an isomorphism. 

It admits several explicit  constructions. 

\subsubsection{First construction: blowing up}\label{sec:first}
 $S_5$ is the blow up of $\BP^2$ at $4$ points $p_1, \ldots, p_4,$ no three of them being collinear. 
$S_5$ contains $10$ straight lines $\ell_i, \ i = 1, \ldots, 10$: $4$ exceptional divisors and the $6$ strict tranforms of the lines 
passing through $p_i, p_j$. 
Their incidence graph, i.e., the graph  with $10$ vertices $v_i$, with $v_i$ being connected with $v_j$ iff 
$\ell_i$ meets $\ell_j$, is the {\it Petersen graph}, \cite{P,D1,D3}. It has $15$ edges. 

\subsubsection{Second construction: Pl\"ucker} 
Consider the Pl\"ucker embedding
$$
Pl: G(2,5) \lra \BP^9 = \BP^{\binom{5}{2} - 1}
$$
We denote the homogeneous coordinates in $\BP^9$ by $x_{ij}, 1\leq i < j \leq 5$. 

Its image is defined by five Pl\"ucker equations
\begin{equation}
x_{ij}x_{kl} - x_{ik}x_{jl} + x_{il}x_{jk} = 0\label{eq:plucker}
\end{equation}
for all $1\leq i < j < k < l\leq 5$. So we have $10$ Pl\"ucker coordinates.

\medskip\noindent{\it Warning:} $\dim G(2,5) = 2\cdot 3 = 6$, so among the $5$ equations \eqref{eq:plucker} only three  are independent.

\begin{claim}\label{claim:section} Let $L\subset \BP^9$ be a $4$-dimensional linear subspace such that 
$$
S:= L\cap Pl(G(2,5))
$$ 
is smooth 
(such $L$ form a Zarisky dense open subset in the Grassmanian of all linear subspaces of dimension $4$). 
Then $S$ is the del Pezzo surface $S_5$.
\end{claim}

See \cite[Prop. 8.5.1]{D3} or \cite[Lemma 2.29]{CKS}. 

\subsubsection{Third construction: the space of stable curves $\overline{M_{0,5}}$} 
Let $M_{0,5}$ be the moduli 
space of curves of genus $0$ with $5$ distinct ordered marked points, i.e. 
$$
M_{0,5} = ((\BP^1)^5 \setminus \bigcup \text{diagonals})/PGL(2).
$$
Let $p_1, \ldots, p_4\in \BP^2$ be in general position, and  
$$
S_5^o = S_5\setminus \bigcup_{i=1}^{10}\ \ell_i\subset S_5
$$
where $\ell_i$ are ten straight lines from section \ref{sec:first}. 

Given $p\in S_5^o$, 
there is a unique conic $C$ passing through $p, p_1, \ldots, p_4$. This conic is isomorphic to 
$\BP^1$ whereupon $5$ distinct points are marked; after passing to the quotient by $PGL(2)$, 
we get a point on $M_{0,5}$.

We get a morphism
$$
S_5^o \lra M_{0,5}
$$
which is an isomorphism. It 
extends to the compactifications, 
$$
S_5 \isom \overline{M_{0,5}}
$$
where $\overline{M_{0,5}}$ is the moduli space of stable curves of genus $0$ with 
$5$ marked points, see \cite{D1}, \cite{K}.

\section{Poncelet-Jacobi curve}
 
\subsection{An elliptic curve related to a couple of plane conics} Let $C, D\subset \BP^2$ be two conics in general position. They meet at four points $p_1, \ldots, p_4$. They give rise to the dual conics 
$C^*, D^*\subset \BP^{2*}$ in the dual projective space, intersection at $4$ points $p_1^*, \ldots p^*_4$, 
corresponding to 4 bitangents of $C, D$.

The Poncelet-Jacobi elliptic curve  is
$$
E = E(C, D) = \{(p, \ell)| \ p\in C,\ p\in \ell, \ell\ \text{tangent to}\ D\}\subset C\times D^* \subset \BP^2\times 
\BP^{2*},
$$
cf. \cite{GH1}. The two projections
$$
C\lla C\times D^*\lra D^*
$$
induce maps

\be \label{eq:CD}
C\overset{\pi_1}\lla E \overset{\pi_2}\lra D^*
\ee
which are both  double coverings ramified at $4$ points, whence we get two involutions 
$i_1, i_2: E\iso E$ whose composition $i_1i_2$ is  a translation by an element $a\in E$. Here 
we have to choose an origin of the group law of $E$ to be 
$$
(p_i, \ell_i)
$$
for some $i$, cf. \cite[\S 2]{GH1}. 

If $a$ happens to be of finite order $n$, the Poncelet process closes up at the $n$-th step.

\subsection{Equation} 

\

An explicit equation for $E$ (due to Cayley) is described in \cite{GH1}. Namely, 
denote the homogeneous coordinates in $\BP^2$ by $x_1, x_2, x_3$. Let two conics $C_1, C_2$ be given by the equations 
$$
\sum_{i,j = 1}^3 c_{ij}x_ix_j = 0,\ \sum_{i,j = 1}^3 d_{ij}x_ix_j = 0.
$$
where $C = (c_{ij}), D = (d_{ij})$ are symmetric $3 \times 3$ matrices. Then $\{ C + \lambda D\}, \lambda\in \BC$ is a pencil of conics passing through $\{ p_1, \ldots, p_4\}$.

Then an equation for $E$ is
\begin{equation}\label{eq:mu}
\mu^2 = \det (C + \lambda D).
\end{equation}

To the regular pentagon there corresponds a degenerate curve $E$. 
If we write a generic curve in the Legendre form 
$$
y^2 = (1 - x^2)(1 - k^2x^2),
$$
then for the degenerate case $k = 0$.

\subsection{Pencil of conics, and another elliptic curve}\label{conics-on-delpezzo} 
We have $5$ forgetful maps
$$
f_i:\ \overline{M_{0,5}} \lra \overline{M_{0,4}} = \BP^1, i = 1, \ldots, 5.
$$
Each $f_i$ has $3$ singular fibers, and the nonsingular fibers 
$$
f_i^{-1}(x)\subset \overline{M_{0,5}} = S_5
$$ 
are conics, and they  form a pencil of conics. 

Three singular fibers correspond to three possible decompositions of the set $\{1, 2, 3, 4\}$ into two pairs of two-element subsets. 

Let $p\in S^o_5$; through $p$ passes the conic $C(p)$ of this pencil. It intersects 
$4$ of the lines $\ell_i$. 

Let $c_1, \ldots c_4$ be the points of intersection of $C(p)$ with these four lines. 
The two-sheeted covering of $C(p)$ ramified at $c_i$ is an elliptic curve, to be denoted $E'(p)$. 

\section{Vinberg's $K3$ mirificum}

\subsection{Definition}\label{vinberg-def}

The Vinberg's "most algebraic $K3$ surface"\ $X_4$ has been introduced and studied in depth in \cite{V}. 
$X_4$ is the Vinberg's notion; we will use the notation $V = X_4$ as well.

It may be defined as follows, cf. \cite[2.1]{DBGKKW}, \cite[5.1]{GT}. 

We start with the union of $6$ straight lines in $\BP^2$ with homogeneous coordinates $(x:y:z)$ given by the equation 
$$
xyz(x - y)(x - z)(y - z) = 0; 
$$
we consider the double covering $\overline P$ of $\BP^2$ branched along these six lines 
(remark that one meets this object in the theory of KZ equation). 

These lines meet three at a time at four points 
$$
p_1 = (1:0:0), p_2 = (0:1:0), p_3 = (0:0:1), p_4 = (1:1:1), 
$$
and our surface $X_4$ is obtained by blowing up $\overline{P}$ at these points.

\subsection{Relation to $S_5$}  On the other hand we can act in the reverse order, so to say, cf. [DBGKKW]. 
 
We blow up $\BP^2$ at $4$ points to get $S_5$. Then we consider the double cover  
$$
\pi_1: \bar S_5 \lra S_5
$$
ramified along its ten straight lines $\ell_i$. The divisor 
$$
\sum_{i=1}^{10} \ell_i
$$
represents the class $-2K_{S_5}$ in $Pic(S_5)$. 

The ten lines $\ell_i$ intersect at $15$ points, and the surface $\bar{S_5}$ is singular at their inverse images, let us denote them $\bar{p_i}, i = 1, \ldots, 15$. 
 
If we blow all $\bar p_i$ we get a smooth surface 
$$
\pi_2:\ \tilde{S_5} \lra \bar{S_5}
$$
which is isomorphic to $X_4$. 

\subsection{Analog of Petersen graph on $X_4$}\label{vinbergs-definition} Originally $X_4$ was defined 
as the unique $K3$ surface whose lattice of transcendental cycles has the form
$$
T = \left(\begin{matrix} 2 & 0\\
0 & 2
\end{matrix}\right),
$$
see \cite{V}, 2.1. Let 
$$
S = Pic(X_4)\subset H_2(X_4, \Bbb Z) 
$$
be the lattice of algebraic cycles; it has the Minkowski signature $(1, 19)$, so 
$$
O(S) = O_{19,1}(\Bbb Z).
$$
The structure of this group has been investigated in \cite{VK}.
Let $O_r(S) \subset O(S)$ be the subgroup generated by reflections with respect to hyperplanes.  
It admits $25$ generators $s_1, \ldots, s_{25}$, and relations are described by the corresponding Coxeter scheme. 

This Coxeter scheme can be seen on \cite{V}, Fig. 1, or on \cite{VK}, Fig. 2 (a). We see that 
this scheme is a refinement of the Petersen graph, so we may suspect that $X_4$ is somehow related 
to $S_5$, and this is so indeed, as we have seen above.

It is not difficult to see $25$ lines in $X_4$ which correspond to reflections $s_i$. Namely, they are 
$10$ strict transforms of the lines $\ell_i\subset S_5$ with respect to the projection 
$$
\pi = \pi_1\pi_2: X_4 \lra S_5,
$$
together with $15$ exceptional divisors.

\subsection{$V$ versus $S_5$} Cf. \cite[2.1]{DBGKKW}. As in 3.4 let $p_1, \ldots, p_4\in \BP^2$ be four points in general position.
The del Pezzo $S_5$ is their blowing up:
$$
b:\ S_5 \lra \BP_2.
$$
Let  $\{ C_t\}$ be the pencil of conics in $\BP^2$ passing through $p_1, \ldots, p_4$, 
and let $\tilde{C}_t\subset S_5$ be the strict transform of $C_t$; we get a pencil of conics $\{ \tC_t\}$ in $S_5$.

Let
$$
\tilde{p}_{t,i} = b^{-1}(p_i) = \tC_t\cap \ell_i\in \tilde{C}_t\subset S_5,\ i = 1, \ldots 4.
$$

Denote $V = X_{-1}$. We have a commutative square

$$
\begin{tikzcd}
  V \arrow[r, "f_2"] \arrow[d, "b_2",swap] & S'' \arrow[d, "b_1"] \\
  S' \arrow[r, "f_1"]                & S_5
\end{tikzcd}
$$

Here $b_1$ is the blowup of $S_5$ at $15$ points $q_1, \ldots, q_{15}$ of intersection of lines $\ell_i,\ i = 1, \ldots, 10$. 

The map $f_1: S'\lra S_5$ is the double cover of $S_5$ ramified along the divisor 
\newline $D = \sum_{i=1}^{10} \ell_i$. 

Let $\ell'_i\subset S''$ be the strict transform of $\ell_i$. The lines $\ell'_i$ are disjoint. 

The map $f_2: V\lra S''$ is the double cover of $S'_5$ ramified along the divisor 
\newline $D' = \sum_{i=1}^{10} \ell'_i$.

The map $b_2$ is the blowup of $S'$ at  points $f_1^{-1}(q_i),\ 1\leq i\leq 15$.

\subsection{A pencil of elliptic curves in $V$}

Let
$$
C'_t := f_1^{-1}(\tilde{C}_t)\subset S'.
$$
It is the double cover of $\tC_t$ ramified at $\tp_{t,i}, i = 1, \ldots, 4$, i.e. an elliptic curve.
Let $E_t\subset V$ be the strict transform of $C'_i$ under $b_1$. 

Alternatively we can go the other way around. Let $C''_t\subset S''$ be the strict transform of 
$\tC'_t$, and $p''_{t,i} = b_1^{-1}(\tp_{t,i})\in C''_t, i = 1, \ldots, 4$. Then 
$$
E_t = f_2^{-1}(C''_t)
$$
appearing now as the double covering of $C''_t$ ramified at $p''_{t,i},  i = 1, \ldots, 4$. 

This way we get map

\be\label{eq:phi}
\phi:\ V\lra \BP^1
\ee

having as non-singular fibers the curves $E_t$, and three singular fibers.

The realisation of $V$ as a double covering shows that $V$ admits an involution 
$i: V\iso V$ which induces the involution on each fiber $E_t$ of the pencil $\phi$. 
 
\begin{claim} The elliptic pencil $\phi$ has three singular fibers of Kodaira type 
${\mathrm I}^*_2$,  or $\tD_6$ in the modern notation.
\end{claim}

In other words, using the notations of \cite[Thm. 6.2]{Ko},  the singular fibers of $\phi$ have the form 
$$
\Theta_0 + \Theta_1 + \Theta_2 + \Theta_3 + \Theta_4 + \Theta_5 + \Theta_6  
$$
where $\Theta_i = \BP^1$, and 
the intersection graph of this divisor has Dynkin type $\tD_6$ (the affine $D_6$), cf. \cite{BHPV}, V, 7. 

\subsection{A combinatorial "Poncelet curve"} 
Consider a regular $n$-gon $P_n$ and a set $E_n$ of pairs 
$(v, e)$ where $v$ is a vertex of $P_n$ and $e$ is an edge of $P_n$ containing $v$. 

The set $E_n$ is a torsor under the dihedral group $D_n$, and we have two projections 
\be\label{eq:proj}
V_n\overset{\pi_{n,1}}\lla E_n \overset{\pi_{n,2}}\lra F_n.
\ee
This diagram is a discrete analogue of the diagram \eqref{eq:CD}.  

\subsection{Two different families of elliptic curves} 
The point $p_i^*\in \BP^{2*}$ may be defined as the variety of straight lines passing through 
$p_i$. 

\

(a)  Let us fix $C$ and vary $D$ in the pencil, so we get a family $(C, D^*_t)$ and corresponding 
family of elliptic curves $E_t\subset C\times D_t^*$. They are all isomorphic since they are 
double coverings of $C$ ramified at $p_1, \ldots, p_4$. 

\

(b)  Let us fix $D$ and vary $C$ in the pencil, so we get a family of elliptic curves 
$\tE_t\subset C_t\times D^*$.  

\

The  curve $\tE_t$ is the double covering of $D^*$ ramified at $p_{1,t}^*, \ldots, p_{4,t}^*$ where 
$p_{i,t}$ are straight lines in $\BP^2$ tangent to $C_t$ and $D$, so $p^*_{i,t}$ varies when $t$ varies.

{\it The curves $\tE_t$ are not isomorphic, in fact they form a universal family.}

The isomorphism class of $\tE_t$ is the cross ratio of $p_{1,t}^*, \ldots, p_{4,t}^*$ in $D^*$.  

\section{Fibonacci pentagrams}

The results of this Section are obtained in collaboration with Gil Bor.

\subsection{Definition of Fibonacci pentagrams}

 A {\em regular} miraculous pentagram has $\phi:=\alpha=\beta=\gamma=\delta=\epsilon$,  so that
$1+\phi=\phi^2$, hence $\phi$ is the {\em golden ratio}, 
$$\phi={1+\sqrt{5}\over 2}.$$

Let 
$$
\phi_n:={F_{n+1}\over F_n},\quad n=1,2,3,\ldots
$$
 be the $n$-th convergent of $\phi$,  where $F_n$ are the Fibonacci numbers, 
 $$
 F_1=F_2=1, F_{n}=F_{n-1}+F_{n-2}.
 $$
 
 Let 2 of the 5 parts of a right triangle be $\phi_n$, say $\gamma_n,\epsilon_n$, then the other 3 parts are determined by equations \eqref{eq:frieze} 
\be
\label{eq:fib-coordinates}
(\alpha_n, \beta_n, \gamma_n, \delta_n, \epsilon_n)=(\phi_{n+1}, \phi_{n+1}, \phi_{n}, {\phi_{n+2}\phi_{n+1}\over \phi_{n}}, \phi_{n}).
\ee
 We call the resulting pentagram  the $n$-th {\em Fibonacci pentagram} $P_n$.
 The  cases $P_1, \ldots,P_4$ are:
\begin{align*}
& (2,2,1,3,1),
\left(\frac{3}{2},\frac{3}{2},2,\frac{5}{4},2\right),
\left(\frac{5}{3},\frac{5}{3},\frac{3}{2},\frac{16}{9},\frac{3}{2}\right),
\left(\frac{8}{5},\frac{8}{5},\frac{5}{3},\frac{39}{25},\frac{5}{3}\right).
\end{align*}

\subsection{Central projection, Poncelet} 

Consider central projection from the sphere onto a  plane tangent to the sphere. It maps great circles to straight lines but  does not preserve angles, in general. There is an exception though. 

\begin{lemma} Under central projection onto a plane tangent to a sphere, the image of two orthogonal great circles, one of which passes through the tangency point, is a pair of orthogonal lines.

\end{lemma}

Now  take a convex self-polar spherical pentagon (or any self-polar  $n$-gon with odd $n$). The diagonals of this pentagon form a smaller pentagon inside it. Fix a point $P$ inside this smaller pentagon and consider the central projection of the spherical polygon onto the tangent plane at $P$. A consequence of the last Lemma is

\begin{cor} The altitudes of the projected planar pentagon are concurrent, intersecting  at $P$.
\end{cor}

Such a pentagon is called {\em orthocentric}.


\subsection{Gauss equation}
The projected planar pentagon  is a Poncelet polygon, i.e. it is inscribed in an ellipse and circumscribes another one. To this pair of ellipses there corresponds an elliptic curve such that the associated  Poncelet map has order 5.

Gauss writes down this curve explicitly. It is written in the Legendre form 
$$
y^2=(1-x^2)(1-k^2x^2)
$$
and will be denoted $E(k)$. The Poncelet map associated with the planar pentagon is given by translation of order 5 on $E(k)$. 

Gauss gives a formula for $k$ in terms of the eigenvalues of the quadratic form in $\BR^3$ whose null cone circumscribes  the given pentagon.  The characteristic equation for the eigenvalues of the quadratic form defining this null-cone is 

\be\label{eq:cubic}
t (2 t-1)^2=(t-1) \omega,
\ee 
where 
$$
\omega=\alpha\beta\gamma\delta\epsilon,
$$
its roots (eigenvalues of the quadratic form) are $G<0<G'< G''$, in terms of which 
\be \label{eq:gauss}
k^2={G'^{-2}-G''^{-2}\over G'^{-2}-G^{-2}}.
\ee

For our Fibonacci polygons $P_n$, 
\be\label{eq:omegan}
\omega_n=\phi_{n}\phi_{n+1}^3 \phi_{n+2}={F_{n+2}^2 F_{n+3}\over F_{n+1}^2 F_{n}}.
\ee

Computer evidence shows that there are simple closed form formulas for the 3 roots, that one of the positive ones is always rational and that $\lim_{n\to\infty}k^2=0.$ 


\subsection{Theorem.}  
(a) $G''$ is rational for $n$ odd, $G'$ is rational for $n$ even. 
(b) It is equal to 
${F_{n+3}/2F_{n+1}}$.

$\Box$

  
The first 10   such rational roots are 
$$
\frac{3}{2},
\frac{5}{4},
\frac{4}{3},
\frac{13}{10},
\frac{21}{16},
\frac{17}{13},
\frac{55}{42},
\frac{89}{68},
\frac{72}{55},
\frac{233}{178}
%
$$

\subsection{•} After dividing equation \eqref{eq:cubic}  by the linear factor given by the rational root, one obtains a quadratic equation for the remaining two roots:
$$
2 t^2 + \frac{ F_n}{F_{n+1}}t-\frac{ F_{n+1}}{F_n}-\frac{ F_n}{F_{n+1}}-2=0,
$$
whose roots are
$$
-\frac{F_n}{4 F_{n+1}}\pm\frac{1}{4} \sqrt{\frac{F_n^2}{F_{n+1}^2}+\frac{8 F_n}{F_{n+1}}+\frac{8 F_{n+1}}{F_n}+16}.
$$

One can now take these explicit expressions for $G,G',G''$ and use equation \eqref{eq:gauss} to get  an explicit formula for $k^2$. 

\newpage

\section{Part II. Fermat's ascent and  Poncelet theorem}

In this Part we compare the Pierre Fermat's ascent method for solving the double equations of Diophantus, as 
described by Andr\'e Weil, with the Poncelet's process. 

\section{The double equations of Diophantus}

The "double equations" appear in the work of Diophantus, and have the form
\begin{equation}\label{dioph-double}
ax^2 + bx + c = u^2,\ a'x^2 + b'x + c' = v^2,
\end{equation}
cf. \cite{Di}, Book IV, Problem 23, \cite{Dif}.

We argue that Pierre Fermat's "ascent method" for solving \eqref{dioph-double}, as 
described by Andre Weil, \cite{W} Chap. II, \S XV, is very close to the Poncelet's construction. 

In both contexts certain elliptic curve plays the key role. For double equations this curve is described in 
\cite{W}, Chap. II. For the Poncelet construction it was discovered by Jacobi, and is nicely described in \cite{GH1}.

\section{Emile Borel's theorem: a rationality criterium}

The following theorem has been used by Bernard Dwork, \cite{D}, in his proof of the rationality of the zeta 
function of algebraic varieties over finite fields.

\begin{theorem}\label{borel-thm} (E. Borel) Let $K$ be a field, 
$$
f(t) = \sum_{i\geq 0} a_it^i\in K[[t]] .
$$
For any $m, n\in \BZ_{\geq 0}$ define an $(m+1)\times (m+1)$ Hankel matrix 
\begin{equation}\label{anm-matrices}
A_{n,m} = \left(\begin{matrix}
a_n & a_{n+1} & \ldots & a_{n+m}\\
a_{n+1} & a_{n+2} & \ldots & a_{n+m+1}\\
. & . & . & \\
a_{n+m} & a_{n+m+1} & \ldots & a_{n+2m}\\
\end{matrix}\right)
\end{equation} 
Let $N_{n,m} = \det A_{n,m}$. Then $f(t)$ is a ratio of two polynomials 
$$
f(t) = \frac{p(t)}{q(t)},\ p(t), q(t)\in K[t],
$$
iff there exist $m, n_0$ such that for all $n\geq n_0$ $N_{n,m} = 0$.
\end{theorem} 

See \cite{B}, \cite{K}, Ch. 5, \S 5, Lemma 5.

\section{Cayley's theorem on the Poncelet porism} 

Consider two smooth conics in $\BP^2$ given by equations $C(x) = 0, D(x) = 0$, $x = (x_0:x_1:x_2)$, 
$C(x) = \sum c_{ij}x_ix_j, \ D(x) = \sum d_{ij}x_ix_j$. 
Consider the formal power series
$$
f(t) = \sqrt{\det(tC + D)} = \sum_{n\geq 0} a_nt^n
$$
the corresponding matrices $A_{n,m}$, \eqref{anm-matrices}, and their determinants $N_{n,m}$ as in the previous Subsection.

\begin{theorem}(A. Cayley) The Poncelet construction yields a finite polygon with $k$ sides if and only if 
$N_{2,m-1} = 0$ for $k = 2m + 1$, or $N_{3,m-2} = 0$ for $k = 2m$. 
\end{theorem}

See \cite{C} or \cite{GH1}.

Let 
\begin{equation}\label{cubic-g}
g(t) = \det(tC + D);
\end{equation}
it is a cubic polynomial.  

Let $D^*\subset \BP^{2*}$ be the dual conic whose points are straight lines $\xi\subset \BP^2$ tangent to $D$. 

The Poncelet - Jacobi elliptic curve $E(C,D)$ is by definition
\begin{equation}
E(C,D) = \{ (x,\xi)\in C\times D^*|\ x\in \xi\}
\end{equation}

\begin{theorem}\label{gh-theorem}
The curve
$$
y^2 = g(t)
$$
is birationally equivalent to $E(C,D)$. 
\end{theorem}

See \cite{GH1}. Another proof will be given below, see \ref{another-pf}.

\section{Intersection of two quadrics in $\BP^3$ and Diophantus double equations}

\subsection{Intersection of two quadrics in $\BP^3$} Following \cite{W}, Chap. II, App. III consider two smooth quadrics $\fC$ and $\fD$ in $\BP^3$ given by equations 
\begin{equation}
\Phi(x) = 0, \Psi(x) = 0,\ x = (x_0:x_1:x_2:x_3).
\end{equation}

We suppose that $\fC$, $\fD$ are in general position;  then $\Omega = \fC\cap \fD$ is a smooth curve. 

Let $\Phi, \Psi$ denote the $4\times 4$ symmetric matrices for $\fC$, $\fD$, and 
$$
F(\xi) = \det(\Phi - \xi\Psi).
$$

\begin{theorem}\label{weil-thm} $\Omega$ is birationally equivalent to the elliptic curve given by 
the equation
$$
\eta^2 = F(\xi).
$$
\end{theorem}

See \cite{W}, Chap. II, App. III.

\subsection{Double equations and Diophantus curve}

Consider a particular case. Let $\fC$ be given by an equation
\begin{equation}\label{double-one}
ax^2 + bxw + cw^2 = u^2,
\end{equation}
and
$\fD$ be given by an equation
\begin{equation}\label{double-two}
a'x^2 + b'xw + c'w^2 = v^2
\end{equation}
in $\BP^3$ with homogeneous coordinates $(w:x:u:v)$. These are Fermat's and Diophantus' "double equations", cf. \cite{W}, 
Chap. II, \S XV, p. 105.

Let $C$ (resp. $D$) be given by the equation \eqref{double-one} (resp. \eqref{double-two}) in $\BP^2$ with coordinates $(w:x:u)$ (resp. $(w:x:v)$); these two guys are isomorphic to $\BP^1$ of course.

We have two maps 
$$
C\overset\alpha\lla \Omega \overset\beta\lra D
$$
where $\alpha(w:x:u:v) = (w:x:u)$ and $\beta(w:x:u:v) = (w:x:v)$,
which make $\Omega$ a double covering of $C$ (resp. $D$) branched at four points. 

Let $i: \Omega\iso \Omega$ (resp. $j: \Omega\iso \Omega$ be the involution which interchanges 
the two branches of $\alpha$ (resp. of $\beta$). Let 
$$
p =  i\circ j:\ \Omega\iso \Omega,
$$
this is an analog of the Poncelet map. 

\subsection{Diophatnus curve is equal to Poncelet curve}\label{another-pf} Let us apply Theorem \ref{weil-thm} to the particular case of double equations. We remark that 
$F(t)$ will coincide with $g(t)$ from \eqref{cubic-g}.

\

We see that the Diophantus curve $E_{Diph} = \Omega$ and the Poncelet curve 
$E_{Ponc}$, cf. \cite{GH1}, p. 40, are given by the same equation (up to replacement $\xi\mapsto  - \xi$). 

Now let us have a look at Fermat's "ascent" method for solving double equations, as described in \cite{W}, Chap. II, \S XV, pp. 105, 106, and 
on an example due to Euler in \cite{W}, Chap. II, App. V, p. 155 - 156. 

Euler discusses an equivalent problem of finding solutions to "canonical equations", equations of the form 
$$
F(\xi, \zeta) = 0
$$
where $F$ is a polynomial of degree $2$ in two variables $\xi, \zeta$. Such an equation defines 
an elliptic curve $C\subset \BP^2$. Euler was completely aware of the relation of this problem to elliptic integrals, cf. \cite{Jlat}.

Starting with some "trivial" solutions $P_0, P_1\in C$, and finding 
a solution $R_0\in C$, he obtains a sequence $R_n\in C$ such that 
$$
R_n \sim R_0 + n(P_1 - P_0)
$$

where for two divisors $D_1, D_2$ on $C$ the notation $D_1\sim D_2$ means that there exists a rational 
function $f$ on $C$ such that 
$$
\text{div}(f) =  D_1 - D_2.
$$ 
The process delivers a finite number of solutions iff $P_1 - P_0$ is of finite order. 

We see that this procedure is similar to the Poncelet process.

\

\section{Exercise} Write down the double equations for a miraculous pentagram.

\end{document}